\documentclass[twoside]{article}

\usepackage{amsthm,amsmath,amssymb,epsfig,verbatim}

{\newtheorem{theorem}{Theorem}[section]

   \newtheorem{proposition}[theorem]{Proposition}

   \newtheorem*{corollary}{Corollary}}

   \newtheorem{propdef}[theorem]{Proposition-Definition}

{\theoremstyle{definition}

		\newtheorem{parag}[theorem]{}

    \newtheorem{problem}[theorem]{Problem}

   \newtheorem*{exercise}{Exercise}

   \newtheorem*{examples}{Examples}}

{\theoremstyle{remark}

  \newtheorem*{definition}{Definition}

    \newtheorem*{remarks}{Remarks}}

\def\>{\mkern1mu}
\def\<{\mkern-1mu}
\def\set{\!:=}

\def\Sing{\text{\rm Sing}}

\def\dt#1{\operatorname{dt_Z}#1}
\def\o#1{{\overline{#1}}}

\def\an{{\text{\rm an}}}
\def\e{\epsilon}
\def\C{\mathbb C}
\def\notC{\raise.25ex\hbox{$\not{}$} \mkern-4.5mu C}
\def\O{\mathcal O}
\def\P{{\mathfrak P}}
\def\p{{\mathfrak p}}
\def\Sigma{{\mathfrak S}}

\title{\bf Equisingularity and Simultaneous Resolution of Singularities}

\author{Joseph Lipman
\thanks{Partially
             supported by the National Security Agency.}}

\date{}

\begin{document}

\maketitle


 \begin{abstract}


Zariski defined equisingularity on an $n$-dimensional
hypersurface~$V$ via stratification by ``dimensionality type,"
an integer associated to a point by means of a generic local projection to
affine \hbox{$n$-space}. A possibly more intuitive concept of equisingularity
can be based on stratification by simultaneous resolvability of
singularities.  The two approaches are known to be equivalent for families of
plane curve singularities. In higher dimension we ask whether constancy of
dimensionality type along a smooth subvariety~$W$ of~$V$ implies the
existence of a simultaneous resolution of the singularities of~$V$
along~$W\<$.  (The converse is false.)

The underlying idea is to follow the classical inductive strategy of
Jung---begin by desingularizing the discriminant of a generic
projection---to reduce to asking if there is a canonical
resolution process which when applied~to quasi-ordinary singularities
depends only on their characteristic monomials. This appears to be so
in dimension~2. In higher dimensions the question is quite open.

\end{abstract}


    \section{Introduction---equisingular stratifications}


The term \emph{equisingularity} has various connotations. Common to
these is the idea of stratifying an algebraic or analytic
$\C\>$-variety~$V$ in such a way that along each stratum the points are,
as singularities of~$V\<$, equivalent in some pleasing sense, and
somehow get worse as one passes from a stratum to its boundary.
A stratification of~$V$ is among
other things a partition into a locally finite family of submanifolds, the
strata (see \S2), so that whatever ``equivalence'' is taken to mean, there
should be, locally on~$V\<$, only finitely many equivalence classes of
singularities.

The purpose, mainly expository and speculative, of this paper---an
outgrowth of a survey lecture at the September 1997 Obergurgl working week---is
to indicate some (not all) of the efforts that have been made to interpret
equisingularity, and   connections among them;  and to suggest directions for
further exploration.


\begin{examples}

(a) The surface~$V$ in~$\C^{\>3}$ given by the equation $X^3=TY^2$ has
singular locus $L\colon X=Y=0$ (the $T$-axis), along which it has
multiplicity~2 except at the origin~$O$, where it has multiplicity 3.
So there is something quite special about~$O$.  Interpreting $V$ along
$L$ as a family of curve germs in the $(X,Y)$-plane, with
parameter~$T$, makes it plausible that the equisingular strata ought
to be $V\!-L$, $L-\{O\}$, and $\{O\}$.
Changing coordinates changes the family, but not the generic
member, and not the fact that something special happens at the
origin.

\smallskip

(b) Replacing $X^3$ by $X^2\<$, we get multiplicity 2 everywhere
along~$L$, including~$O$. The corresponding family of plane curve
germs consists of a pair of intersecting lines degenerating to a
double line; so there is still something markedly special about $O$---a
feeling reinforced by the following picture.

\epsfxsize=2.9in
\bigskip\smallskip
\centerline{\epsfbox{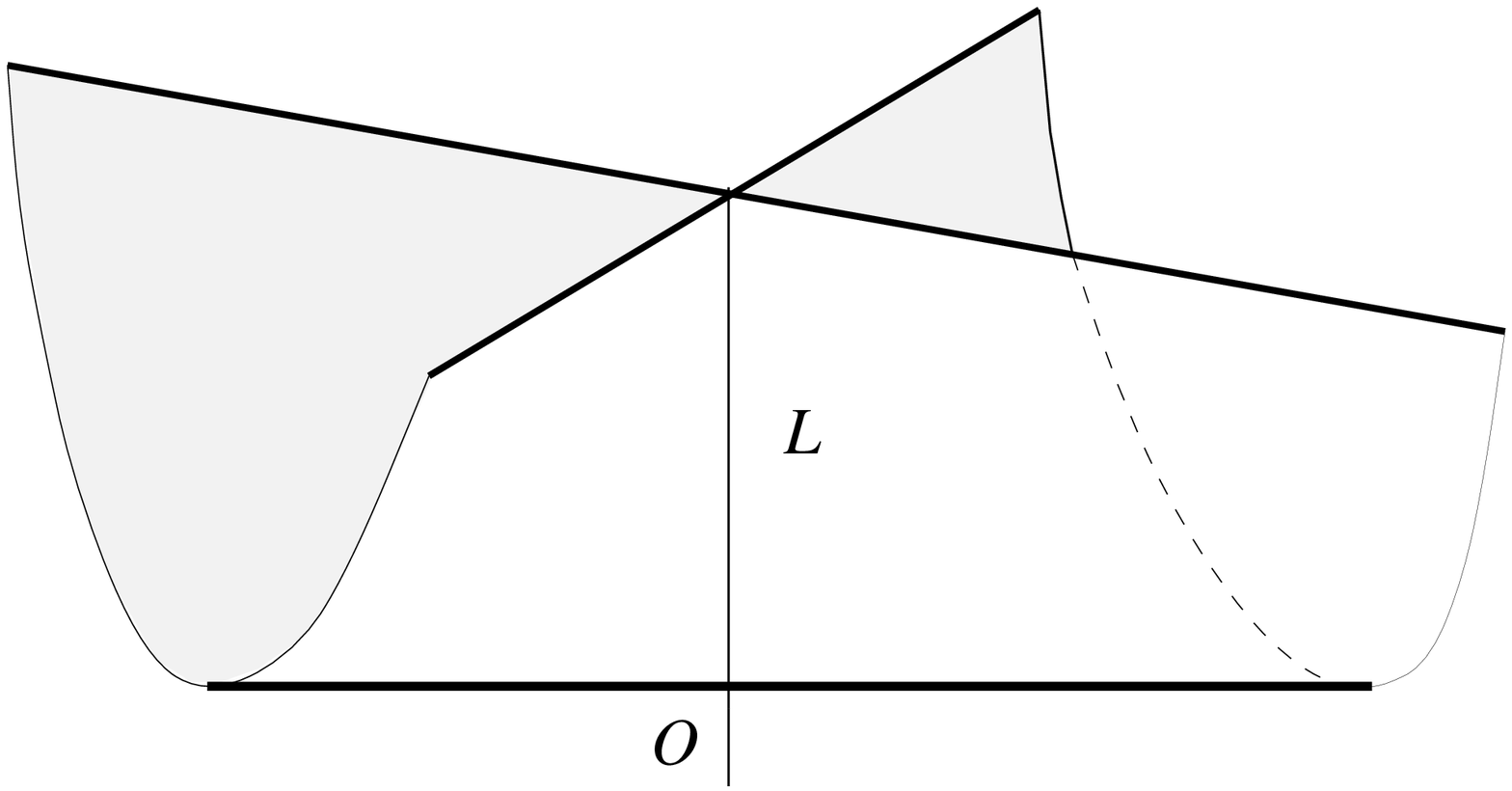}}
\bigskip

This indicates  that equimultiplicity along~$L$, while
presumably \emph{necessary}, is \emph{not sufficient} for
equisingularity.

\smallskip

(c) The surface $V\colon X(X+Y)(X-Y)(TX+Y)=0$, along its singular locus
$L\colon  X=Y=0$, can be regarded as a family of plane-curve germs all of
which look topologically the same, namely four distinct lines through
a point. Thus we have \emph{topological equisingularity} along $L$.
(In this example the homeomorphisms involved can even be made
\emph{bi-Lipschitz;} see
\cite {P} for more on Lipschitz equisingularity.) From this point of view, the
equisingular strata should be
$V\!-L$ and $L$.  However, from the analytic---or even differentiable---point
of view, all these germs are distinct, since the cross-ratio of the four lines
varies with~$T$. Thus differential isomorphism is too stringent a
condition for equisingularity---there are too many equivalence
classes.\looseness=-1

\smallskip

(d) \emph{Differential equisingularity} of $\>V$ along $L$
at a smooth point~$x$ of~$L$ is associated with the Whitney
conditions $\mathfrak W(V\<, \>L)$ holding at~$x$.  These
conditions signify that if $y\in L$ and
$z\in V\!-\Sing(V)$  (where $\Sing(V)$ is
the singular locus of~$V$) approach~$x$ in such a way that the tangent space
$T_{V\!,\>z}$ and the line joining $y$ to $z$ both have limiting positions,
then the limit of the tangent spaces contains the limit of the lines.
(The conditions can be described intrinsically, i.e., as a condition
on the prime ideal of~$L$ in the local ring $\O_{V\!\<,\>x}$, see proof of
Theorem~\ref{Z is W} below.)

\epsfxsize=2.6in
\centerline{\epsfbox{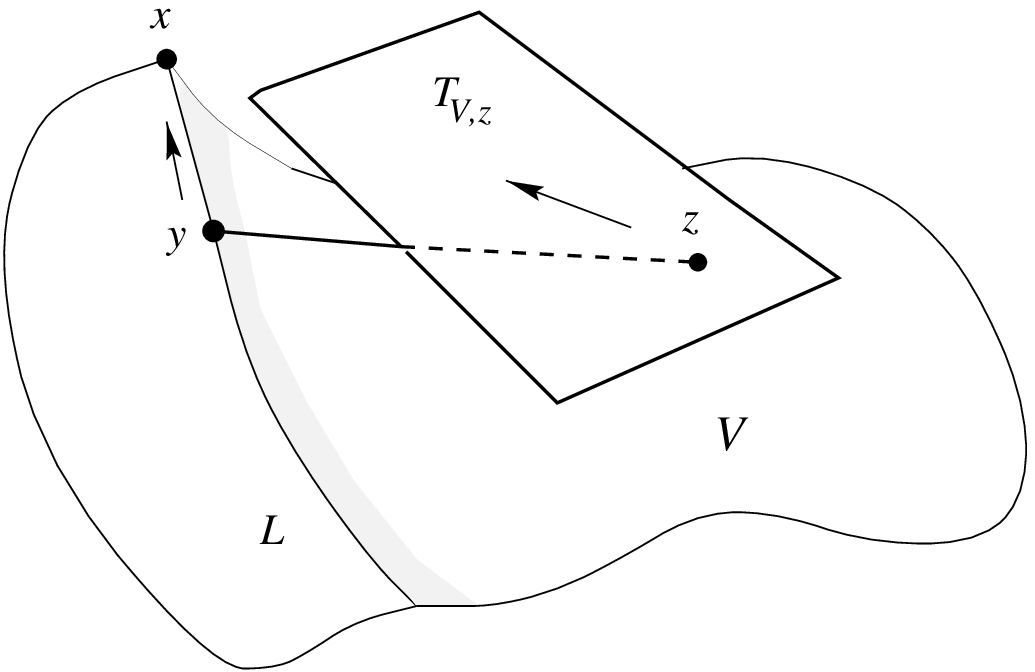}}
\medskip\smallskip

The conditions fail in example (b) above for
$y=(X_y,Y_y,T_y)\!=(0,0,t^{1/2})$ and
$z=(X_z,Y_z,T_z)\!:=(t^2\<,t,t^2)$, the limiting tangent space as
$t\to 0$ being the plane
\hbox{$T=2X$,} while the limiting line is just $L\colon X=Y=0$.

The conditions $\mathfrak W(V\<, \>L)$ by themselves are not enough:
in \cite[p.\,488]{Z1}, Zariski gives the example $X^2=T^2Y^3$ in which
$\mathfrak W(V,L_i)$ holds for $L_1\colon X=T=0$ and $L_2:X=Y=0$, but
topological equisingularity does not hold at the origin along either
of these components of the singular locus. One must require that a
neighborhood of~$x$ in~$L$ be contained in a stratum of some Whitney
stratification---a stratification such that $\mathfrak W(\o{S_1},S_2)$
holds everywhere along~$S_2$ for any strata $S_1$, $S_2$ with
$\o{S_1}\supset S_2$. In fact Thom and Mather showed that a Whitney
stratification of~$V$ is topologically equisingular in the following
sense: for any point~$x$ of a stratum~$S$, the germ $(V,x)$ together
with its induced stratification is topologically the product of
$(S,x)$ with a stratified germ $(V_0,x)\subset (V,x)$ in which $x$ is a
stratum by itself (see \cite[p.\,202, (2.7), and p.\,220,
Corollary~(8.4)]{M1}).

The obvious formulation of a converse to the Thom-Mather theorem is
false: in~\cite{BS}, Brian\c con and Speder showed that the family of
surface germs
$$
Z^5+TY^6Z+Y^7X+X^{15}=0
$$
(each member, for small~$T\<$, having an isolated singularity at the
origin) is topologically equisingular, but not differentially so. On
the other hand, there is a beautiful converse, due to L\^e and
Teissier, equating topological and differential equisingularity of the
totality of members of the family~$(V\cap H,L,x)$ as $H$ ranges over
general linear spaces containing~$L$ (see
\cite[\S5]{LT}, \cite[p.\,480, \S4]{T2}) .

\end{examples}

The approach to equisingularity via $\mathfrak W$ is the most
extensively explored one. Whitney stratifications are of basic
importance in the classification theory of differentiable maps
(\cite{M1}, \cite{M2}, \cite{GL}). (See also  \cite{DM}, for the existence of
Whitney stratifications in a general class of geometric categories.)
 Hironaka proved that a Whitney stratification is \emph{equimultiple:} for any
two strata $S_1$,
$S_2$, the closure~$\o {S_1}$ has the same multiplicity (possibly 0) at every
point of~$S_2$
\cite[p.\,137, 6.2]{H1}.  More generally, work of L\^e and Teissier
in the early 1980's, generalized in the 1990's by Gaffney and others,
has led to characterizations of Whitney equisingularity by the
constancy of a finite sequence of ``polar multiplicities." (For this,
and much else, see, e.g., \cite{T2}, \cite{GM}, \cite{K}).

\penalty -2000
It is therefore hard to envision an acceptable definition of
equisingularity which does not at least \emph{imply} differential
equisingularity.  But one may still wish to think about conditions
which reflect the analytic---not just the differential---structure
of~$V\<$.  Indeed,
\S\S3--5 of this paper are devoted to describing two plausible formulations
of ``analytic equisingularity" and to trying to compare them with each
other and with differential equisingularity.
\smallskip

In brief, we think of equisingularity theory as being \emph{the study
of conditions which give rise to a satisfying notion of ``equisingular
stratification,\kern-1pt" and of connections among them.}


          \section{Stratifying conditions}


To introduce more precision into the preceding vague indications about
equisingularity (see e.g.,~\ref{wrapup}), we spend a few pages on
basic remarks about stratifications. Knowledgeable readers may prefer
going directly to \S3.

There is a good summary of the origins of stratification theory
in \hbox{pp.\:33--44} of \cite{GMc}.
What follows consists mostly of
variants of material in \cite[pp.\:382--402]{T2},
and is straightforward to
verify. The main result, Proposition~\ref{coarsest}, generalizes
\cite[pp.\:478--480]{T2} (which treats the case of Whitney stratifications).


\begin{parag} \label{parts}

We work either in the category of reduced complex-analytic spaces or
in its subcategory of reduced finite-type algebraic $\mathbb
C$-schemes; in either case we refer to objects~$V$ simply as
``varieties."  For any such~$V$ the set $\Sing(V)$ of singular
(non-smooth) points is a closed subvariety of~$V\<$.  A
\emph{locally closed subvariety of\/~$V$} is a subset~$W$ each of
whose points has an open neighborhood $U\subset V$ whose intersection
with $W$ is a closed subvariety of~$U\<$. Such a $W$ is open in its
closure.

By a \emph{partition} (analytic, resp.~algebraic) of a variety~$V$ we
mean a locally finite family $(P_\alpha)$ of non-empty subsets
of~$V$ such that for each~$\alpha$ both the closure~$\o{P_\alpha}$ and
the boundary $\partial P_\alpha\set\o{P_\alpha}-P_\alpha$ are closed
subvarieties of~$V$ (so that $P_\alpha$ is a locally closed subvariety
of~$V$), such that $P_{\alpha_1}\cap P_{\alpha_2}=\emptyset$ whenever
$\alpha_1\ne\alpha_2\>$, and such that $V=\cup_\alpha \>P_\alpha\>$. The
$P_\alpha$ will be referred to as the \emph{parts} of the partition.

For example, a partition may consist of the fibers of an
upper-semicontinuous function $\mu$ from $V$ into a well-ordered
set~$I$, where ``upper-semicontinuous" means that for all $i\in I$,
$\{\,x\in V\mid\mu(x)\ge i\,\}$ is a closed subvariety of~$V\<$.

We say that a partition $\P_1$ \emph{refines} a partition $\P_2\>$,
or write $\P_1\!\prec\P_2\>$, if, $\P_1$ and~$\P_2$ being identified
with equivalence relations on~$V\<$---subsets of~$V\<\times V\<$---we
have \hbox{$\P_1\subset \P_2$.}  So
$\P_1\!\prec\P_2\>$ iff the following (equivalent) conditions
 hold:\vadjust{\kern1.5pt}

(i) Every $\P_1$-part is contained in a $\P_2$-part.\vadjust{\kern1.5pt}

(i)$'$ With $\P_x$ denoting the $\P$-part containing~$x$,
$\P_{1\<,x}\subset\P_{2\<,x}$ for all $x\in V\<$.\vadjust{\kern1.5pt}

(ii) Every $\P_2$-part is a union of~$\P_1$-parts.

\penalty -1000

Let $\P$ be a partition of~$V\<$, and let $W\subset V$ be an
\emph{irreducible} locally closed subvariety. For any $\P$-part $P\<$,
$P\>\cap W=\o P\cap W\!-\partial P\cap W$ is the difference of two
closed subvarieties of~$W\<$, so its closure in~$W$ is a closed
subvariety of~$W\<$, either equal to $W$ or of lower dimension
than~$W$; and $W$ is the union of the locally finite family of all
such closures, hence equal to one of them.  Consequently there is a
unique part---denoted $\P_W$---whose intersection with~$W$ is dense
in~$W\<$. Moreover, $W-\P_W=\partial\P_W\cap W$ is a proper closed
subvariety of~$W\<$.

\medskip
\end {parag}

\begin{parag} \label{strats}


Condition (ii) in the following definition of stratification is
non-standard---but suits a discussion of equisingularity (and
may be necessary for Proposition~\ref{coarsest}).

\begin{definition}
A \emph{stratification} (analytic, resp.~algebraic) of~$V$ is a
decomposition into a locally finite disjoint union of non-empty,
connected, locally closed subvarieties, the \emph{strata,} satisfying:
\begin{enumerate}
\item[(i)] For any stratum~$S$, with closure~$\o S$, the boundary $\partial
S\set\o S-S$ is a union of strata (and hence is a closed subvariety
of~$V$).
\item[(ii)] For any $\o S$ as in (i),  $\Sing(\overline S)$ is a union
of~strata.
\end{enumerate}
\end{definition}

{}From (i) it follows that a stratification is a partition, whose parts
are the strata.  Noting that a subset $Z\subset V$ is a union of
strata iff $Z$ contains every stratum which it meets, and that
$\Sing(\overline S)$ is nowhere dense in~$\overline S$, we deduce:
 \begin{enumerate}
\item[(iii)] Every stratum is smooth.
\end{enumerate}

The strata of codimension~$i$ are called $i$-strata.

\smallbreak
For stratifications $\Sigma_{\<1}$, $\Sigma_2$, we have
$\Sigma_{\<1}\!\prec\Sigma_2\Leftrightarrow$ the closure of any
$\Sigma_2$-stratum is the closure of a $\Sigma_{\<1}$-stratum.

\medskip

\end {parag}


\begin{parag} \label{filts}

 One way to give a stratification is via a filtration
$$
V=V_0\supset V_1\supset V_2\supset\dots
$$
by closed subvarieties such that:

\begin{enumerate}
\item[(1)] For all $i\ge 0$, $V_{i+1}$ contains $\Sing(V_i)$, but contains no
component of~$V_i\>$.  In other words, $V_i-V_{i+1}$ is a dense
submanifold of~$V_i\>$.
\end{enumerate}
The strata are the connected components of $V_i-V_{i+1}\ (i\ge0)$, and
their closures are the irreducible components of the~$V_i\>$. (Note
that (1) forces the germs of the~$V_i$ at any $x\in V$ to have
strictly decreasing dimension, whence $\cap\,V_i$ is empty.) It
follows that a closed subset $Z$ of $V$ is a union of strata
$\Leftrightarrow$ for~any component~$W$ of any~$V_i\>$, if
$W\not\subset Z$ then $W\cap Z\subset V_{i+1}\>$.  So for (i) and (ii)
above to hold we need:

\begin{enumerate}
\item[(2)]
If $W\<$, $W'$ are irreducible components of~$V_i\>$, $V_j$
respectively and $W\not\subset W'$ then $W\cap W'\subset V_{i+1}\>$,
and if $ W\not\subset \Sing(W')$ then $W\cap \Sing(W')\subset
V_{i+1}\>$.
\end{enumerate}

For any stratification, let $0=n_0<n_1<n_2<\dots$ be the integers
occurring as codimensions of strata. Redefine $V_i\set$ union of all
strata of codimension $\ge n_i\>$, thereby obtaining a filtration
which gives back the stratification and which satisfies:

\begin{enumerate}
\item[(3)] There is a sequence of integers $0=n_0<n_1<n_2<\dots$ such that
$V_i$ has pure codimension $n_i$ in~$V\<$.
\end{enumerate}

\noindent Thus there is a one-one correspondence between the set of
stratifications and the set of filtrations which satisfy (1), (2) and (3).

\medskip

In the following illustration, a tent which should be imagined to be
bottomless and also to stretch out infinitely in both horizontal
directions, $V_1$ can be taken to~ be the union of the ridges (labeled
$E$, $F$), and $V_2$ to be the vertex $O$.

\epsfxsize=3.0 in
\bigskip\smallskip
\centerline{\epsfbox{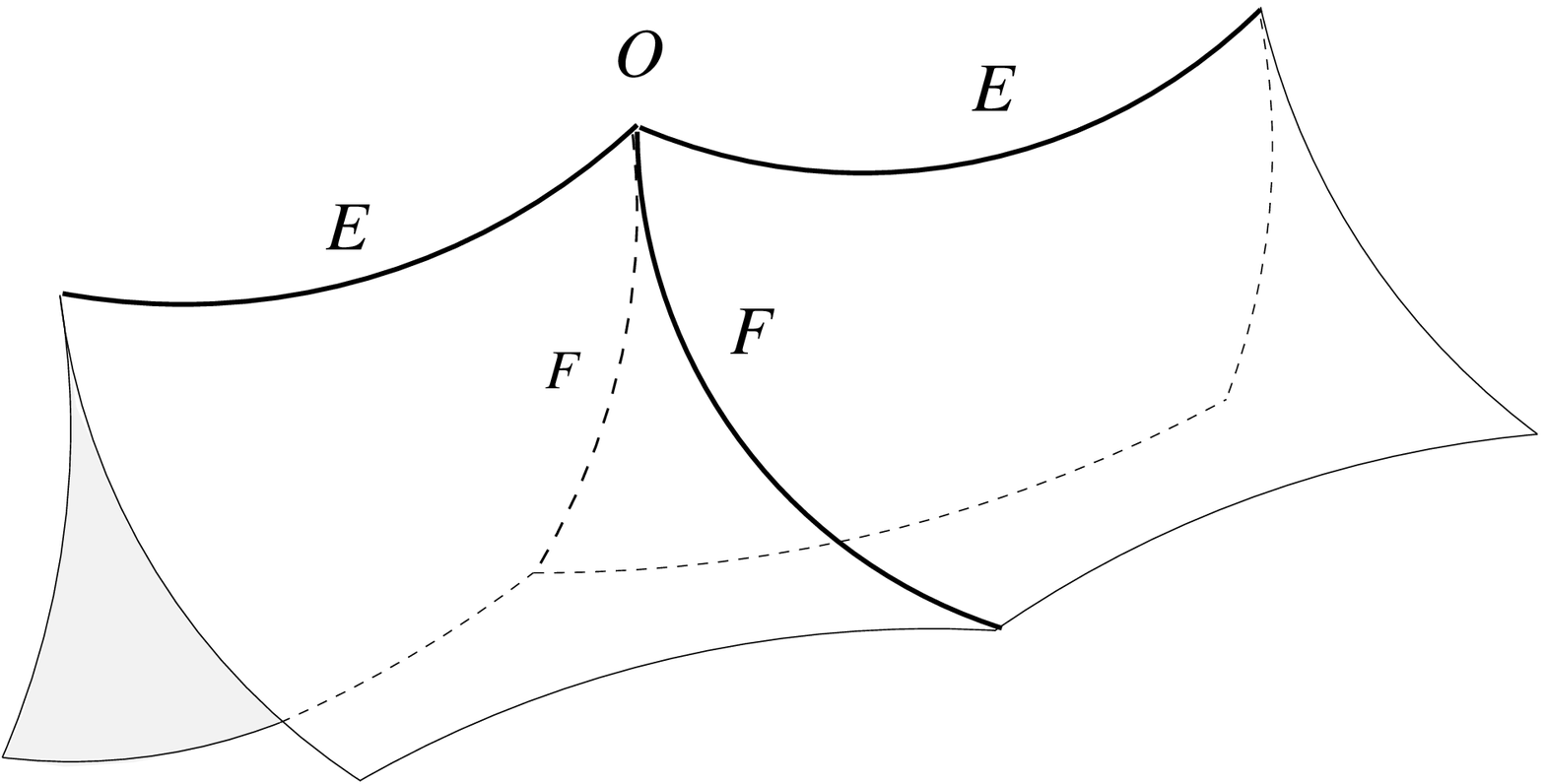}}
\bigskip

  To say that an irreducible subvariety $W$ is not contained in a
closed subvariety~$Z$ is to say that $W\cap Z$ is nowhere dense
in~$W\<$. It follows that the conditions (1), (2) and~(3) are local:
they hold for a filtration if and only if they hold for the germ of
that filtration at every point~$x\in V$---and indeed, they can be
checked inside the complete local ring~$\widehat{\mathcal
O_{V\!,\>x}}$. Moreover if a filtration $\mathfrak F_{\<x}$ of a germ
$(V,x)$ satisfies these conditions, then there is such a
filtration~$\mathfrak F_{\<N}$ of an open neighborhood~$N$ of~$x$
in~$V\<$, whose germ at~$x$ equals $\mathfrak F_{\<x}$. Thus there is
a \emph{sheaf of stratifications} whose stalk at~$x$ can be identified
with the set of all filtrations~$\mathfrak F$ of~$(V,x)$ satisfying
(1), (2) and~(3).\looseness=-1
\medskip
\end{parag}


\begin{parag} \label{scondn}


Stratifications can be determined by \emph{local stratifying conditions,} as
follows. We~consider conditions $C=C(W_1,W_2,x)$ defined for all $x\in
V$ and all pairs \mbox{$(W_1,x)\supset (W_2,x)$} of subgerms\-
of~$(V,x)$ with $(W_1,x)$
\emph{equidimensional}  (all components
of~$W_1$ containing~$x$ have the same dimension) and $(W_2,x)$
\emph{smooth}.  Such a pair can~be thought of as two radical~ideals
$p_1\subset p_2$ in the local ring $R\set\O_{V\!\<,\>x}\>$, with
$R/p_1$ equidimensional and $R/p_2$ regular. For~example,
$C(W_1, W_2, x)$ might be defined to mean that $W_1$ is equimultiple
along~$W_2 $ at~$x\>$---i.e., the local rings $R/p_1$ and
$(R/p_1)_{p_2}$ have the same multiplicity.  Or, $C(W_1, W_2, x)$
could signify that the Whitney conditions $\mathfrak W(W_1,W_2)$ hold
at~$x$, these conditions being expressible within~$\O_{V\!\<,\>x}$.\looseness=1

(Recall that there is a natural equivalence between the category of
local analytic \hbox{$\mathbb C$-algebras} and the category of complex
analytic germs---the category of pointed spaces localized with respect
to open immersions, i.e., enlarged by the adjunction of formal
inverses for all open immersions. A similar equivalence holds in the
algebraic context.)

For such a~$C$ and for any subvarieties~$W_1\>$, $W_2$ of~$V$ with
$W_1$ closed and locally equidimensional, and $W_2$ locally closed,
set
\begin{align*}
C(W_1,W_2) & \set
\bigl\{\,x\in W_2\,\bigl|\, W_2
\text{ smooth at $x$,  and if $x\in W_1$ then} \\
 &\smash{\phantom{\set \{\,x\in W_2\,\bigl|\ \ }
\text { $\bigl(\>(W_1,x)\supset (W_2,x)$ and
$C(W_1,W_2 ,x)\>\bigr)$}\,\bigr\}. }\\
\notC(W_1,W_2)&\set  W_2-C(W_1,W_2). {\vphantom{\Bigl(}}
\end{align*}

\vspace{-3pt}
\noindent For example, if $W_1$ contains no component of~$W_2$ then
$$
\notC(W_1,W_2)=\Sing(W_2)\cup (W_1\cap W_2).
$$
The condition~$C$ is called \emph{stratifying} if for any such $W_1$ and
$W_2\>$, $\notC(W_1,W_2)$ is contained in a nowhere dense closed
subvariety of~$W_2\>$. (It suffices that this be so whenever $W_2$ is
smooth, connected, and contained in~$W_1$.)

For any $S\subset V\<$, denote by $S^\an$ the smallest closed
subvariety of~$V$ which contains~$S$.  If $S$ is a union of
$\Sigma$-strata for some stratification $\Sigma$, then $S^\an$ is just
the topological closure of~$S$.

A stratification~$\Sigma$\vadjust{\kern.5pt}  is a
\emph{$C$-stratification} if for any strata~$S_1$, $S_2$ with
$S_2\subset\o{S_1}$ it holds that
$\notC(\>\o {S_1},\o{S_2}\>)^\an$ is a union of strata.  For
\emph{any}\vadjust{\kern.5pt} two  strata $S_1$, $S_2$ of a $C$-stratification,
$\notC(\>\o {S_1},\o{S_2}\>)^\an$ is a union of strata, since
$$
\o{S_1}\not\supset S_2\implies
\notC(\>\o{S_1},\o{S_2}\>)=\Sing(\o{S_2})\cup(\>\o {S_1}\cap\o{S_2}\>);
$$
and if moreover $C$~is stratifying then $C(\>\o {S_1},S_2)=S_2\>$.
For some stratifying~$C$, this last condition by itself may or may
not be enough to make $\Sigma$ a
$C$-stratification; if it always is enough, then we'll say that
$C$ is a \emph{good} stratifying condition.
(See, for instance, Example (d) below.)
If $C$ is good and if $C'$ is any stratifying
condition which implies~$C$, then every $C'$ stratification is a
$C$-stratification.

\smallskip
We will show in~\ref{defsigma} below that for each stratifying
condition~$C$ there exists a $C$-stratification. In fact there exists
a \emph{coarsest} $C$-stratification---one which is refined by all
others (see~\ref{coarsest}).

\end{parag}


\begin{parag} \label{xampl}


\begin{examples}

(a) If $C$ is the empty condition (i.e., $C(W_1,W_2,x)$ holds for all
pairs $(W_1,x)\supset (W_2,x)$ as in \S\ref{scondn}) then $C$ is a good
stratifying condition, and every stratification is a $C$-stratification.

\smallskip

(b) The logical conjunction (``and") of a finite family of stratifying
conditions, \smash{$\hat C\set\wedge_{i=1}^n\, C_i\>$},
 is again stratifying.
If a stratification~$\Sigma$ is a $C_i$-stratification for all~$i$
then $\Sigma$ is a $\hat C$-stratification; but the converse does not
always hold.  It does hold if each~$C_i$ is good---in which case $\hat C$
is good too.

\smallskip
\penalty -1500

(c) Let $\P=(P_\alpha)$ be a partition of~$V\<$.  For  $x\in V$
let $\alpha_x$ be such~that $x\in P_{\alpha_x}.$ Define
$C_\P(W_1,W_2, x)$ to mean that $(W_2,x)\subset (P_{\alpha_x},x)$. Then
$C_\P$ is stratifying: as in~\ref{parts}, for any
irreducible locally closed subvariety $W_2\subset V$ there is a
unique~$\alpha$ such that $W_2-P_\alpha$ is a proper closed
subvariety of~$W_2$, and then for any closed locally
equidimensional $W_1\supset W_2$ we have
$\notC_\P(W_1,W_2)=(W_2-P_\alpha)\cup\Sing(W_2)$,
a nowhere dense closed subvariety of~$W_2\>$.
The condition $C_\P\<(\>\o{S_1},S_2)$ on a pair of strata of a
stratification~$\Sigma$ means that the part $\P_{S_2}$ of~$\P$ which
meets $S_2$ in a dense subset actually contains~$S_2$,
and hence is the only part of~$\P$ meeting~$S_2\>$.
One deduces that \emph{a $C_\P$-stratification is just a stratification
which refines~$\P$.}\looseness=-1


\begin{exercise}

For two partitions $\P$, $\P'\<$, $\>\bigl(C_{\P'}\Rightarrow
C_\P\bigr)\iff \P'\!\prec\P$.

\end{exercise}


(d) More generally, suppose given a partition~$\P_W$ of~$W$ for every
equidimensional locally closed subvariety~$W\subset V\<$. (For example, if
$(P_\alpha)$  is a partition of~$V\<$, one can set $\P_W\set (W\cap
P_\alpha)$.) With $C_W\set C_{\P_W}$ as in Example~(c),
define $C(W_1, W_2,x)$ to mean $C_{W_1}(W_1,W_2,x)$. Arguing
as before, one finds that $C$ is a stratifying condition
on~$V$, and that a
$C$-stratification is one such that, if $S$ is any stratum then each
part of~\smash{$\,\P_{\o{\!S\<}}\,$} is a union of strata---a
condition which holds if and only if for any two strata $S_1$, $S_2$, we have
$C(\>\o{S_1},S_2)=S_2$. Thus this~$C$ is
a good stratifying condition.

\smallskip

(e) The above-mentioned Whitney conditions $\mathfrak W(W_1,W_2,x)$
are stratifying. This was first shown, of course, by Whitney,
\cite[p.\,540, Lemma 19.3]{Wh}. It was shown much later, by Teissier,
that ${\raise.25ex\hbox{$\not{}$}\mkern-8mu {\mathfrak W}}(W_1,W_2)$
itself is analytic (\cite[p.\,477, Prop.\,2.1]{T2}.  Indeed, the
main result in \cite{T2}, Theorem~1.2 on p.\,455, states in part that
the stratifying condition $\mathfrak W$ is of the type described in
example~(d) above, the partition $\P_W$ being given by the
level sets of the polar multiplicity sequence on~$W\<$.

\end{examples}


\end{parag}


\begin{propdef} \label{defsigma}

For any stratifying condition\/ $C,$\ the following inductively
defined filtration of\/ $V$ gives rise, as in~\ref{filts}, to a\/
$C$-stratification~$\Sigma_C\colon $ $ V_0=V,$ and for\/ $i>0,$\
$V_{i+1}$ is the union of all the $W^*\subset V_i$ such that

\text{\rm (a)} $W^*$ is a component either of\/ $\Sing(V_j)$
for some\/ $j\le i$ or of\/ $\notC(W'\<,W)^\an$ for some components\/
$W'\<,$ $W$ of\/~$V_j\>,$ $V_k$ respectively, with $j\le k\le i\>;$\
and

\text{\rm (b)} $W^*$ is not a component of\/ $V_i\>$.

\end{propdef}


\begin{proof}

It is clear that $V_{i+1}$ contains every component of $\Sing(V_i)$,
but no component of~$V_i\>$. If $W$ and $W'$ are components of~$V_i$
and $V_j$ respectively, and if $j>i$, then $W\cap W'\subset V_{i+1}$;
while if $j\le i$ and $W\not\subset W'$ then
\hbox{$W\cap W'\subset\notC(W'\<, W)\subset V_{i+1}$,} since by the
definition of stratifying condition, no component of~$\notC(W'\<,
W)^\an$ is a component of~$V_i\>$. So, strata being connected
components of~$V_i-V_{i+1}\ (i\ge0)$, the closure of any
stratum---i.e., any component of any~$V_i\>$---is a union of strata
(see~\ref{filts}).  Now any component~$W^*$ of $\Sing(V_j)$ is a
component of~$V_i$ for some~$i>j\>$: otherwise, by the definition
of~$V_i\>$, $W^*\subset\cap V_i=\emptyset$ (see~\ref{filts}). Hence
$\Sing(V_j)$ is a union of strata.  Similarly, if $W'\supset W$ are
closures of strata, and so components of~$V_j\>,$ $V_k$ respectively,
with $j\le k$, then $\Sing(W'\<,W)^\an$ is a union of strata. Thus the
filtration does indeed give a $C$-stratification.
\end{proof}

\penalty -1500

\begin{remarks}

1. If $C(W,W,x)$ holds for every smooth point $x$ of every
component~$W$ of~$V\<$, then $V_1 =\Sing(V)$.

2. Any stratification~$\Sigma$ is $\Sigma_C$ for a good~$C$, viz.\
 $C\set C_\Sigma$ (Example~\ref{xampl}(c)).  This results at once from
 Proposition~\ref{coarsest} below.  Or, if $\mathfrak V\colon
 V=V_0^*\supset V_1^*\supset\dots$ is the filtration associated
 to~$\Sigma$, satisfying conditions~(1), (2) and~(3) in~\ref{filts},
 so that
\hbox{$V_i^*\subset V$} has pure codimension, say~$n_i\>$, then
one can check for irreducible components $W\<$, $W'$ of~$V_i^*\>$,
$V_j^*$ respectively, that
$$
\notC(W'\<,W)=\notC(W'\<,W)^\an\subset W\cap V_{i+1}^*=\notC(W,W),
$$
whence a component of $\notC(W'\<,W)$ is a component of~$V_\ell^*$ iff
it is the closure of an $n_\ell$-stratum on the boundary of~$W$; and
it follows that the filtration described in
Proposition-Definition~\ref{defsigma} is identical with $\mathfrak
V\<$.

\end{remarks}


\begin{proposition} \label{coarsest}

For any stratifying condition\/ $C,$\ $\Sigma_C$ $($defined
in~$\ref{defsigma})$ is the coarsest\/ $C$-stratification
of\/~$V\<$---every\/ $C$-stratification refines\/~$\Sigma_C\>$. In particular,
if\/ $C$ is\/
\emph{good} \textup{(\S\ref{scondn})} then $\Sigma_C$ is the coarsest
among all stratifications such that\/ $C(\o{S_1},S_2)=S_2$ for any
two strata\/ $S_1,$ $S_2$.

\end{proposition}


\begin{proof}

Let $V=V_0\supset V_1\supset V_2\supset\dots$ be as in
{}~\ref{defsigma}, and let $\Sigma$ be a $C$-stratification of~$V\<$.
The assertion to be proved is: \emph{If\/ $Z$ is an irreducible
component of\/ $V_j$ $(j\ge 0),$\ then\/ $Z=\o{\Sigma_{\<\<Z}},$\ the
closure of\/ $\Sigma_{\<\<Z}$.}  (Recall from~\ref{parts} that $\Sigma_{\<\<Z}$
is
the unique stratum containing a dense open subset of $Z$, and see the
last assertion in~\ref{strats}.)

For $j=0$, it is clear that $Z=\o{\Sigma_{\<\<Z}}\>$.  Assume the assertion
for all $j\le i$.  Let $Z$ be a component of~$V_{i+1}\>$.
By~\ref{defsigma}, $Z$ is a component either of~(a): $\Sing(V_j)\
(j\le i)$ or of (b): $\notC(W'\<,W)^\an$ with $W'$ a component
of~$V_j$ and $W$ a component of~$V_k$ $(j\le k\le i)$.  In case~(a),
the inductive hypothesis gives that every component~$W''$ of~$V_j$ is
the closure of a $\Sigma$-stratum, so both $W''$ and $\Sing(W'')$ are
unions of $\Sigma$-strata (see~\ref{strats}); and it follows easily
that $\Sing(V_j)$ is a union\vadjust{\kern.4pt} of $\Sigma$-strata, one of
which
must be~$\Sigma_{\<\<Z}$ (because $\Sing(V_j)$ meets $\Sigma_{\<\<Z}$).
Since\vadjust{\kern.7pt}
$ Z\subset\o{\Sigma_{\<\<Z}}\subset\Sing(V_j) $ and $Z$~is a component
of~$\Sing(V_j)$, therefore $Z=\o{\Sigma_{\<\<Z}}$.  Similarly,
in\vadjust{\kern.4pt} case~(b) the inductive hypothesis gives that $W'$ and $W$
are both closures of
$\Sigma$-strata, and so, $\Sigma$~being a $C$-stratification,
$\notC(W'\<,W)^\an$ is a union of $\Sigma$-strata, one of which must
be $\Sigma_{\<\<Z}$, whence, as before $Z=\o{\Sigma_{\<\<Z}}$. Thus the
statement
holds for $j=i+1$, and the Proposition results, by induction.
\end{proof}


\begin{corollary}

 If\/ $C$ and\/ $C'$ are stratifying conditions such that\/ $C$ is
good and\/ \mbox{$C'\Rightarrow C$}  then\/ $\Sigma_{C'}\<\<\prec\Sigma_C\>$.

\end{corollary}



\begin{corollary}

{\rm (Cf.~\cite[p.\,536, Thm.\,18.11, and p.\,540, Thm.\,19.2]{Wh}.)}
If\/ $C$ is a good stratifying condition and\/
$\>\P_1,\P_2,\dots,\P_n$  are  partitions of\/~$V\<,$ there is a coarsest\/
$C$-stratification refining all the\/ $\P_i$.

\end{corollary}


\begin{proof}

The stratifications~$C_{\P_{\<i}}$ of Example~\ref{xampl}(c)\vadjust{\kern.3pt}
are good (see Example~\ref{xampl}(d)), and so in view of
Example~\ref{xampl}(b)),\vadjust{\kern.3pt} Proposition~\ref{coarsest} shows
that $\,\Sigma_{C\wedge C_{\P_{\<\<1}}\wedge\dots\wedge C_{\P_{\<\<n}}}$
does the job.
\end{proof}

\penalty -1000

\begin{parag} \label{wrapup}


 \kern-1.5ptFor a  stratifying condition~$C$ say that $V$ is
\emph{$C$-equisingular} at a subgerm~$(W,x\<)$ if there is a
$C$-stratification~$\Sigma$ such that $x\in
\Sigma_W$, the unique $\Sigma$-stratum containing a dense open subvariety
of~$W\<$---or equivalently, if $(W,x)\subset(\Sigma_{C,x\>},x)$ where
$\Sigma_{C,x}$ is the unique $\Sigma_C$-stratum containing~$x$.

As a special case of the exercise in Example~\ref{xampl}(c), it holds
for any~two stratifying conditions $C$ and $C'$ on\/~$V$ that
\begin{align*}
\bigl(\Sigma_{C'}\<\<\prec\Sigma_C\bigr)\iff&
\bigl(\text{for all subgerms $(W,x)$
of~$V\<$,}\\ V\< &\text{ is $C'\<$-equisingular at~$(W,x)\Rightarrow
V\<$ is $C$-equisingular at~$(W,x)$}\bigr).
\end{align*}
In particular, by the first corollary of Proposition~\ref{coarsest}, if $C$ is
good
and
$C'\Rightarrow C$ then
$C'\<$-equisingularity implies $C$-equisingularity.

\end{parag}


\section {The Zariski stratification}

In the early 1960's, Zariski developed a comprehensive theory of
equisingularity in codimension~1 (see \cite{Z3}). Let $x$ be a point
on a hypersurface~$V\<$, around which the singular locus~$W$ is a
smooth manifold of codimension~1.  The fibers of any local retraction
of~$V$ onto~$W$ form a family of plane curve germs. Zariski showed
that if the singularity type of the members of one such family is
constant---where singularity type is determined in the classical sense
via characteristic pairs, embedded topology, multiplicity sequence,
etc.---then the
\emph{Whitney conditions} hold along $W$ at $x$; and conversely, the Whitney
conditions imply \emph{constancy of singularity type} in any
such family.  In this case, moreover, the Whitney conditions are
equivalent to \emph{topological triviality} of~$V$ along~$W$
near~$x$.  Furthermore, if these conditions hold, the singularities
of~$V$ along $W$ near $x$ can be resolved by blowing up~$W$ and its
successive strict transforms (along which the conditions continue to
hold), in a way corresponding exactly to the desingularization of any
of the above plane curve germs by successive blowing up of infinitely
near points---and conversely. (This last condition can be thought of
as \emph{simultaneous desingularization} of any family of plane curve
germs arising as the fibers of a retraction.)

Unfortunately, in higher codimensions no two of these characterizations
of equi\-singularity in codimension~1 remain equivalent; and anyway
there is as yet no definitive sense in which two hypersurface
singularities can be said to have the same singularity type.
Exploration of the remaining connections among these characterizations
leads to interesting open questions, a few of which are stated
below.\looseness=-1

After some tentative attempts at generalizing  the notion
of equi\-singularity to higher codimensions (see e.g., \cite[p.\,487,
\S4]{Z1}), Zariski formulated (in essence) the following definition of
the dimensionality type $\dt(V,x)$ of a hypersurface
germ~$(V,x)$---possibly empty---of dimension~$d$, where $d=-1$ if
$(V,x)$ is empty, and otherwise $x$ is a $\C\>$-rational point of the
algebraic or analytic variety~$V\<$. (One can regard such a $(V,x)$
algebraically as being a local $\C\>$-algebra of the form $R/(f)$ where
$R$ is a formal power-series ring in $d+1$ variables over~$\C$ and
$0\ne f\in R$.)  Zariski-equisingularity is defined by local constancy
of~$\dt{}$.

The intuition which inspired Zariski's definition does not lie on the
surface.  It may have come out of his extensive work on ramification
of algebraic functions and on fundamental groups of complements of
projective hypersurfaces, or have been partly inspired by Jung's method of
desingularization (which begins with a desingularization---assumed, through
induction on dimension, to exist---of the discriminant of a general
projection.)  Anyway, here it is:
\begin{align*}
&\dt(V,x)=-1&&\text{if $(V,x)$ is empty,}\\
&\dt(V,x)=1+\dt(\Delta_\pi,0)\qquad&&\text{otherwise,}
\end{align*}
where $\pi\colon (V,x)\to (\C^{\>d}\<,0)$ is a general finite map germ,
with discriminant~$(\Delta_\pi,0)$ (the hypersurface subgerm of
$(\C^{\>d}\<,0)$ consisting of points over which the fibers of~$\pi$ have
less than maximal cardinality, i.e., at whose inverse image $\pi$ is
not
\'etale). Here $\pi$ is defined by its coordinate functions
$\xi_1,\dots,\xi_d$, power series in $d+1$ variables, with linearly
independent linear terms; and if $(V,x)$ is then represented---via
Weierstrass preparation---by an equation
$$
f(Z)=Z^n+a_1(\xi_1,\dots,\xi_d)Z^{n-1}+\cdots+a_n(\xi_1,\dots,\xi_d)=0
\qquad (a_i\in\C\mkern.6mu[[\xi_1,\dots,\xi_d]])
$$
then $\Delta_\pi$ is given by the vanishing of the $Z$-discriminant
of~$f\<$, so that $(\Delta_\pi,0)$ is a hypersurface germ of dimension
{}~$d-1$, whose $\dt{}$ may be assumed, by induction on dimension,
already to have been defined.  A property of finite map germs $\pi$
holds for
\emph{almost all}~$\pi$ if there is a finite set of polynomials in the
(infinitely many) coefficients of $\xi_1,\dots,\xi_{d\>}$ such that the
property holds for all~$\pi$ for whose coefficients these polynomials
do not simultaneously vanish.  (The coefficients depend on a choice of
generators for the maximal ideal of the local $\C\>$-algebra
of~$(V,x)$---i.e., of an embedding of~$(V,x)$ into~$\C^{\>d+1}\<$---but
the notion of ``holding for almost all~$\pi\>$" does not.) It follows
from \cite[p.\,490, Proposition~5.3]{Z2} that $\dt(\Delta_\pi,0)$ has
the same value for almost all~$\pi\>$; and that enables the preceding
inductive definition of $\dt{}$ (which is a variant of the original
definition in~\cite[\S4]{Z2}).

It is readily seen that $(V,x)$ is smooth iff $\Delta_\pi$ is empty
for almost all~$\pi$, i.e., $\dt(V,x)=0$.  So $\dt(V,x)=1$ means that
almost all $\pi$ have smooth discriminant, which amounts to Zariski's
classical definition of codimension-1 equisingularity of $(V,x)$ (see
e.g., \cite[pp.\:20--21]{Z3}.

For another example, suppose the components~$(V_i,0)$
 of~$(V,0)\subset(\C^{\>d+1}\<,0)$ are smooth, with \emph{distinct}
 tangent hyperplanes\vadjust{\kern.7pt} \mbox{$\sum_{j=1}^{d+1}
 a_{ij}X_j=0$}.  Then a similar property holds for
 $(\Delta_\pi,0)=(\cup_{i\ne i'}\>\>\pi(V_i\cap V_{i'}),0)$ if $\pi$
 is general enough; and it follows by induction on~$d\>$ that $\dt(V,0)$
 is one less than the rank of the matrix~$(a_{ij})$.

\smallskip

Hironaka proved, in \cite{H2}, that on a $d$-dimensional
\emph{algebraic} $\>\C$-variety $V$ which is locally embeddable
in~$\C^{\>d+1}$, \emph{the function~$\dt{}$ is upper-semicontinuous.}
More precisely, it follows from the main Theorem on p.\,417
of~\cite{H2} and from
\cite[p.\,476, Proposition~4.2]{Z2} that the set
$$
\postdisplaypenalty 10000
V_i\!=\{\,x\in V\mid \dt(V,x)\ge i\,\}\qquad(i\ge0)
$$
is a closed subvariety of~$V\<$. Thus we have a partition~$\P_{\<\<Z}$
of~$V\<$, by the connected components of the fibers of the
function~$\dt{}$, and correspondingly a stratifying condition
$\mathfrak Z(W_1,W_2,x)$, defined to mean that $\dt{}$ is constant on
a neighborhood of~$x$ in~$W_2\>$. This is a  good stratifying condition,
see Example~\ref{xampl}(d). We call the stratification
$\Sigma_{\<\mathfrak Z}$ the \emph{Zariski stratification}. By
Proposition~\ref{coarsest}, the Zariski stratification is the coarsest
stratification with $\dt{}$ constant along each stratum.

One would expect to have a similar stratification for
\emph{complex-analytic} locally-hypersurface varieties; but this is
not explicitly contained in the papers of Zariski and Hironaka. To be
sure, I asked Hironaka during the workshop (September, 1997) whether
his key semi-continuity proof applies in the analytic case, and he
said not necessarily, there are obstructions to overcome. Thus:

\begin{problem}
Investigate the upper-semicontinuity of\/ $\dt{}$ on analytic
hypersurfaces, and the possibility of extending the Zariski
stratification to this context.
\end{problem}

Zariski proved in \cite[\S6]{Z2} that with $V_i$ as above,
$V_i-V_{i+1}$ is smooth, of pure codimension~$i$ in~$V\<$
\cite[p.\,508, Theorem~6.4]{Z2}; and that the closure of a part
of~$\P_{\<\<Z}$ is a union of parts [\emph{ibid.,} pp.\:510--511].  One
naturally asks then whether $\P_{\<\<Z}=\Sigma_{\<\<Z}$---what is still
missing is condition~(ii) of Definition~\ref{strats}, that the singular locus
of
the closure of a part is a union of parts.  But that is indeed so, as
follows from the equimultiplicity assertion in the next Theorem.


\begin{theorem} \label{Z is W}
Let\/ $V$ be a purely\/ $d$-dimensional algebraic\/ $\C\>$-variety,
everywhere of embedding dimension\/ $\le d+1$.
Then for any two parts\/ $P_1,$\ $P_2$ of the
partition\/~$\P_{\<\<Z}$\vadjust{\kern.5pt} with
$\o{P_1}\supset P_2,$\ the Whitney conditions\/~$\mathfrak W(\o{P_1},P_2,x)$
hold at all\/~$x\in P_2\>$. So $\o{P_1}$ is
equimultiple along~$P_2,$ and\/
$\P_{\<\<Z}=\Sigma_{\<\<Z}$ is a Whitney stratification of\/ $V\<$.
\end{theorem}


\begin{remarks}

(a) Let us say that  $V$ is
\emph{Zariski-equisingular} along a subvariety~$W$ at a smooth
point~$x$ of~$W$ if $\dt{}$ is constant on a neighborhood of~$x$
in~$W$ (see \cite[p.\,472]{Z2}). The equality $\P_{\<\<Z}=\Sigma_{\<\<Z}$
entails
that Zariski-equisingularity is the same as $\mathfrak
Z$-equisingularity (see~\S\ref{wrapup}).

(b) From Theorem~\ref{Z is W} we see  via Thom-Mather, that
Zariski-equisingularity implies topological triviality of~$\>V$ along $W$
at~$x$,
a result originally due to Varchenko (who actually needed only \emph{one}
sufficiently general  projection with discriminant topologically trivial along
the
image of~$W$), see
\cite{Va}.

(c) In connection with the idea that \emph{analytic} equisingularity
should involve something more than differential equisingularity, note
the example of Brian\c con and Speder
\cite[p.\,3, (2)]{Z3} showing that $\dt{}$ need not be
constant along the strata of a Whitney stratification. In fact, in
\cite{BH} Brian\c con and Henry show that for families of isolated
singularities of surfaces in~$\C^{\>3}\<$, constancy of~$\dt{}$ is equivalent
to
constancy of the generic polar multiplicities (i.e., Whitney equisingularity)
and of two \emph{additional} numerical invariants of such singularities.

\end{remarks}
\goodbreak


Here is a sketch of a \emph{proof of Theorem~\ref{Z is W}}.
We first need a formulation of~$\mathfrak W(W_1,W_2,x)$ ($x$ a smooth point of
$W_1\subset W_2$)
in terms of the local
$\C\>$-algebra $R\set\O_{W_{\<1}\<,\>x}$ and the prime ideal~$\p$ in~$R$
corresponding to~$W_2$ (so that $R/\p$ is a regular local ring).  One fairly
close
to the geometric definition is, in outline, as follows. We assume for
simplicity that $W_1$ is equidimensional, of dimension~$d$.  Let
$m$ be the maximal ideal of~$R$, and in the Zariski tangent space
$T\set \textup{Hom}_{R/m}(m/m^2\<, R/m)$ let $T_2$ be the tangent space
of~$W_2\>$, i.e., the subspace  of  maps in~$T$ vanishing on
$(\p+m^2)/m^2$.  Set $W_1'\set \textup{Spec}(R)$,
$W_2'\set \textup{Spec}(R/\p)$, let $B\to W_1'$ be
the blowup of~$W_1'$ along~$W_2'\>$, let $\Omega$ be
the universal finite differential module of $W_1'/\C$, and
let $G\to W_1'$ be the Grassmannian which functorially
parametrizes rank $d\/$ locally free quotients of~$\Omega$.  Then the
canonical map $p\colon B\times_{W_1'} G\to {W_1'}$ is an
isomorphism over $V_0\set {W_1'}\<-W_2'\<-\Sing({W_1'})$;
and we let
$Y\subset B\times_{W_1'} G$ be the closure of~$p^{-1}V_0$.  Any $\C\>$-rational
point~$y$ of the closed fiber of $Y\to W_1'$ gives rise naturally, via
standard universal properties of~$B$ and $G$, to a pair $(E,F)$ of
vector subspaces of~$T$, where $E$ contains $T_2$ as a codimension-1
subspace, and $\dim F=d\>$: $E$ consists of all maps in~$T$ vanishing
on the kernel of the surjection
$(\p+m^2)/m^2\cong\p/m\p\twoheadrightarrow\mathfrak q$, where
$\mathfrak q$ is the 1-dimensional quotient corresponding to the
$W_1'$-homomorphism
\mbox{$\textup{Spec}(R/m)\to B\subset \mathbb P(\p)$} whose  image is the
projection of~$y$,  $\mathbb P(\p)$ being the projective bundle
parametrizing \hbox{1-dimensional} quotients of~$\p$; and $F$ is dual to a
similarly-obtained $d$-dimensional quotient of the $R/m$-vector space
$\Omega/m\Omega\cong m/m^2$. Then $\mathfrak W(W_1,W_2,x)$ holds
$\iff$\vspace{5pt}

\noindent
 $\mathfrak W(R,\p)$: \emph{$E\subset F$ for all\/
$(E,F)$ arising from points in the closed fiber of\/ $Y\to W_1'$}.

\smallskip

\noindent (The idea, which can be made precise
with some local analytic
geometry, is that for any embedding of $W_1$ into $\>\C^{\>t}\ (t=\dim T)$
with $W_2$ linear, the set of pairs $(E,F)$ can naturally be identified
with the set of limits of sequences $(E_n,F_n)$ defined as follows:
let $x_n\to x$ be a sequence in~$W_1-W_2-\Sing(W_1)$, let $E_n$ be the join
of~$W_2$ and~$x_n$---a linear space containing~$W_2$ as a 1-codimensional
subspace---and let $F_n$ be the tangent space to $W_1$ at~$x_n$.)

It is not hard to see that
$\,\mathfrak W(R,\p)\iff \mathfrak W(\hat R,\p\hat R)$.

This leads to a reduction of the proof of~Theorem~\ref{Z is W} to that
of the similar statement for the \emph{formal} Zariski partition of
$\textup{Spec}(\widehat{\O_{V\!,\>x}})$. (That partition is by constancy of
formal dimensionality type, defined inductively through
\emph{generic} finite map germs---those whose coefficients (see above)
are independent indeterminates.)  Indeed, we have to prove $\mathfrak
W(R,\p)$ whenever $R$ is the local ring of the closure of a
Zariski-stratum~$W_1$ of~$V$ at a $\C\>$-rational point~$x$, and $\p$ is
the prime ideal corresponding to the Zariski-stratum~$W_2$
through~$x$. (We refer---prematurely---to the
parts of the Zariski partition as Zariski-strata.)  The dimensionality type of
the generic point of~$\textup{Spec}(R)$ is equal to that of all the
closed points of~$W_1$, i.e., to the codimension of~$W_1$ in~$V$ (see
\cite[p.\,419, Theorem]{H2} and
\cite[p.\,508, Thm.\,6.4]{Z2}); and the analogous statement follows for the
generic point of each component of~$\textup{Spec}(\hat R)$
(see~\cite[p.\,507, bottom]{Z2})---so that each of those components is
the closure of a formal Zariski stratum
of~$\textup{Spec}(\widehat{\O_{V\!,\>x}})$ (\cite[p.\,480,
Thm.\,5.1]{Z2}).  To complete the reduction, note that if
\mbox{$\widehat{R_1},\dots,\widehat{R_s}$} are the quotients of~$\hat R$ by
its minimal primes then \looseness=-1
$$
\mathfrak W(\hat R,\p\hat R)\iff
\mathfrak W(\widehat{R_i},\p\widehat{R_i})
\quad\textup{for all }i=1,2,\dots,s.
$$
(The corresponding geometric statement is clear.)

Once we are in the formal situation, where $V=\textup{Spec}
(\C[[X_1,\dots,X_{d+1}]]/(f))$ with $f$ an irreducible formal power
series, \cite[p.\,490, Prop.~5.3]{Z2} gives an embedding
of~$V$ into $\C^{\>d+1}$ such that for almost all \emph{linear}
projections \mbox{$\pi\colon\C^{\>d+1}\<\<\to\C^{\>d}\mkern-1.5mu$,} if $w$
is either the generic point of~$W_1$ or the generic point of~$W_2$ or
the closed point of~$V$ then the dimensionality type of the point
$\pi w$ on the discriminant hypersurface~$\Delta_\pi$ of~$\pi|V$ is
one less than that of~$w$, while the dimension of the closure $\o{\pi w}$ is
the
same as the dimension of~$\o w$. It follows then from \cite[p.\,491,
Prop.~5.4]{Z2} that $\pi(\o{W_i})\ (i=1,2)$ is the closure of a Zariski stratum
of~$\Delta_\pi$, except when $\o{W_i}=V\<$.

At this point we have reached a situation which is much like the one
treated by Speder in~\cite{S}. In the context of local analytic
geometry, Speder shows that a certain version of
equisingularity does imply the Whitney conditions. That version is not
the same as the one we have called Zariski-equisingularity: while it
is also based on an induction involving the discriminant, the
induction is with respect to
\emph{one} sufficiently general linear projection,
whereas Zariski's induction is
with respect to \emph{almost all} projections. Nevertheless there are
enough similarities in the two approaches that Speder's arguments can
now be adapted to give us what we want.  Note that Speder asserts only
that equisingularity gives the Whitney conditions when $\o{W_1}=V\<$;
but in view of the remark about~$\pi(\o{W_i})$ in the preceding
paragraph, the remaining cases can be treated \emph{in the present
situation} by reasoning like that in section~IV of~\cite{S}.

Details (copious) are left to the reader.


\section {Equiresolvable stratifications}


A standard way to
analyze singularities is to {\it resolve\/} them, and on the resulting
manifold to study the inverse image of the original singularity. For
example, after three point blowups,
the cusp at the origin of the plane curve $Y^2=X^3$ is replaced by a
configuration of three lines crossing normally, a configuration which
represents the complexity of the original singularity.

\epsfxsize=2.6 in
\bigskip
\centerline{\epsfbox{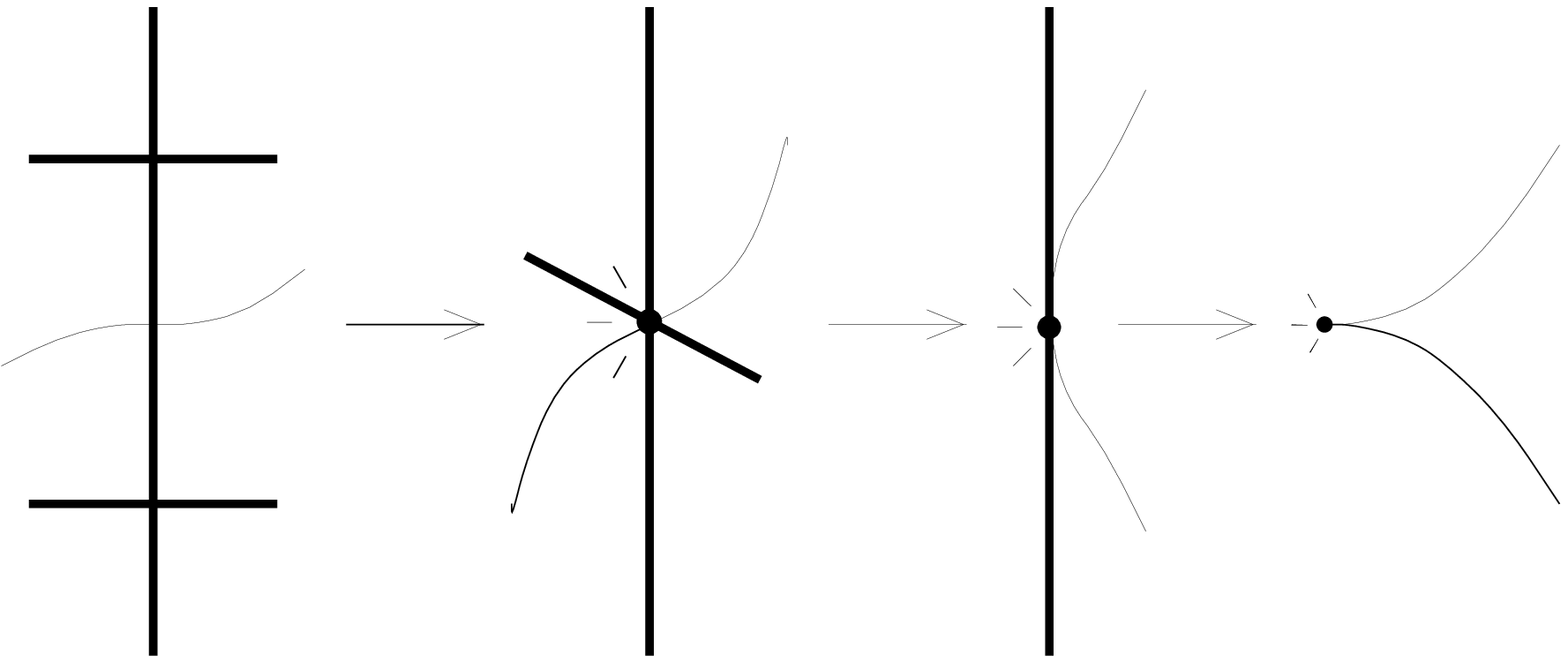}}

\medskip

Quite generally, the  embedded
topological type of any plane plane curve singularity
is determined by the intersection matrix of its inverse image---a
collection of normally crossing projective lines---on a minimal
embedded desingularization.

 This suggests one sense in which a family of singularities could be
regarded as being equisingular: if its members can be resolved
simultaneously,\nopagebreak[4] in such a way that their inverse images all
``look the same."

Several reasonable ways to make such an idea precise come to mind.
Here is one. With notation as in section~\ref{scondn}, let us  define the
\emph{local equiresolvability} condition $ER(W_1,W_2,x)$ to mean that there
exists an embedding of the germ $(W_1,x)$ into a smooth germ~$(M,0)$, a
proper birational (or bimeromorphic) map of
manifolds $f\colon M'\to M$ such that $f^{-1}W_1$ is a divisor with
smooth components having only normal crossing intersections---so that
$f^{-1}W_1\>$ (considered as a \emph{reduced} variety) has an obvious
stratification, by multiplicity---and such that $f^{-1}W_2$ is a union of
strata,
each of which $f$ makes into a $C^\infty\<$-fiber bundle, with smooth fibers,
over~$W_2$.

\begin{proposition}
{\rm (Cf.~\cite[p.\,401]{T2}.)} The condition\/ $ER\/$ is stratifying (see
\S\ref{scondn}).
 Hence, by  Proposition~\ref{coarsest}, every variety has a coarsest
equiresolvable stratification.

\end{proposition}

\begin{proof}

Let  $W_2\subset W_1\subset V\<$, with $W_2$ smooth and
connected. The question being local, we may embed $W_1$ in some open
$M\subset\C^{\>n}\<$. Let $f_2\colon M_2\to M$ be the blowup of
\mbox{$W_2\subset M$,} so that
$f_2^{-1}W_2$ is a divisor in~$M_2\>$. Let $f_1\colon M'\to M_2$ be an
embedded resolution of~$f_2^{-1}W_1$, and set $f\set f_2f_1$,
so that $f^{-1}W_1\subset M'$ is a finite union of normally crossing smooth
divisors, as is  its subvariety $f^{-1}W_2 $ (also a divisor in~$M'$). After
removing some proper closed subvarieties from~$W_2$ we may assume that
each intersection~$W$ of components of~$f^{-1}W_2$ is mapped onto~$W_2$.
Then  the theorem of Bertini-Sard allows us to assume further  (after removing
another proper closed subvariety of~$W_2$) that  the restriction of~$f$ maps
each such~$W$ \emph{submersively} onto~$W_2\>$. Now the first isotopy
lemma of Thom \cite[p.\,41]{GMc}, or \cite[p.\,311, Thm.\, 4.14]{Ve},
applied to the multiplicity stratification---clearly a Whitney
stratification---of~$f^{-1}W_2\>$, shows that each stratum is,
via $f\<$, a smooth $C^\infty\<$-fiber bundle over $W_2\>$.\looseness=-1
\end{proof}


Following Teissier, \cite[p.\,107, D\'efinition\,3.1.5]{T1}, we can define a
stronger
equiresolvability condition~$ER^+$ by adding to~$ER\/$ the requirement that
the \emph{not-necessarily-reduced} space $f^{-1}W_2$ be locally
analytically trivial over each $w\in W_2$, i.e., germwise the product of~$W_2$
with the fiber~$f^{-1}w$. \emph{Condition~$ER^+\/$ is still stratifying:} with
notation as in the preceding proof, set $E_2\set f_2^{-1}W_2\>$, and  realize
the desingularization $f_1$ as a composition of blowups
$f_{i+1}\colon M_{i+1}\to M_i\ \mbox{$(i\ge 2)$}$ of smooth
$C_i\subset M_i$ such that if $E_{j+1}\set f_j^{-1}(C_j\cup E_j)\
(j\ge 2)$, then $C_i$ has normal crossings with  $E_i$; then note for any $z\in
M_i$ at which both $C_i$ and  $(f_2f_3\dots f_i)^{-1}W_2$ are  locally
analytically trivial, that $(f_2f_3\dots f_if_{i+1})^{-1}W_2$ is locally
analytically trivial everywhere along $f_{i+1}^{-1}z$; and finally argue as
before,
using Bertini-Sard etc.

The result that $ER^+$ is stratifying is similar to \cite[p.\,109,
Prop.\,2]{T1}; but  Teissier's formulation of strong
simultaneous resolution refers to birational rather than embedded resolution,
and so seems a priori weaker than~$ER^+\<$. I don't know (but would like to)
whether in fact Teissier's condition is equivalent to~$ER^+\<$. \looseness=-2

I also don't know whether either of the stratifying conditions $ER$ and $ER^+$
is~\emph{good} (see \S\ref{scondn}), or in case they are not, whether there is
a
convincing formulation of equiresolvability which does lead to a good
stratifying
condition.

\penalty -1500

Teissier showed that strong simultaneous resolution along any smooth
subvariety $W$ of the hypersurface~$V$
implies the Whitney conditions\nopagebreak[3] for the pair
$(V-\text{Sing}(V),W)$, see \cite[p.\,111, Prop.\,4]{T1}.

A converse, Whitney $\Rightarrow$ strong simultaneous resolution,
was proved by  Laufer (see~\cite{L}) in case $V$ is the total space of a family
of isolated two-dimensional surface singularities,for example if $V$ is
three-dimensional, $\text{Sing}(V)=:W$ is a nonsingular curve, and there is a
retraction $V\to W$ whose fibers are the members of the family.

\medskip
The basic motivation behind this paper is our interest in possible relations
between Zariski-equisingularity and simultaneous resolution.  Such
relations played a prominent role in Zariski's thinking---see e.g.,
\cite[p.\,490, F, G. H]{Z1} and \cite[p.\,7]{Z3}.

To begin with, Zariski showed that if $\dt(V,x)=1$ then
$\text{Sing}(V)$ is smooth, of codimension~1 at~$x$, and that for
some neighborhood $V'$ of~$x$ in~$V\<$, the Zariski stratification
of~$V'$ \emph{is} equiresolvable (in a strong sense, via successive
blowups of the singular locus, which stays smooth of codimension~1
until it disappears altogether)---see \cite[p.\,93, Thm.\,8.1]{Z3}). Zariski's
result involves birational resolution, but can readily be extended to cover
embedded resolution.

The \emph{converse,} that any simultaneous resolution of a family of plane
curve singularities entails equisingularity, was shown by Abhyankar \cite{A}.
This converse fails in higher dimension,  for instance for
the Brian\c con-Speder family of surface singularities referred to in
Remark~(c) following Theorem~\ref{Z is W}, a family which is differentially
equisingular and hence---by the above-mentioned result of Laufer---strongly
simultaneously resolvable, but not Zariski-equisingular.

We mention in passing a different flavor of work, by Artin, Wahl, and others,
on
simultaneous resolution  and infinitesimal equisingular deformations of normal
surface singularities, see e.g., \cite{Wa}.

\smallskip
At any rate, a positive answer to the next question would surely enhance the
appeal  of defining analytic equisingularity via~$\dt$.


\begin{problem} \label{Z is E}

Is the Zariski stratification of a hypersurface~$V$  equiresolvable in some
reasonable sense? For example, is it an\/
$ER$- or\/~$ER^+\/$-stratification?

\end{problem}


Initial discouragement is generated by an example of Luengo \cite{Lu},  a
family of quasi-homogeneous two-dimensional hypersurface singularities
which is equisingular in Zariski's sense, but cannot be resolved simultaneously
by blowing up smooth centers. There is however another way. That is
the subject of the next section.


\section{Simultaneous resolution of quasi-ordinary singularities}


We pursue the question ``Does Zariski-equisingularity imply
equiresolvability?"  (A positive answer would
settle Problem~\ref{Z is E} if equiresolvability
were expressed by a good stratifying condition, see~\S\ref{wrapup}.)

In fact we ask a little more: if $x\in V\<$,  with
$V$ a  hypersurface  in~$\C^{\>d+1}\<$, does there exist a neighborhood $M$ of
$x$
in~$\C^{\>d+1}$ and a proper birational (or bimeromorphic) map of manifolds
$f\colon M'\to M$ such that
$f^{-1}V$ is a divisor with smooth components having only normal crossing
intersections and such that for \emph{every}
Zariski-stratum~$W$ of~$V\cap M$, $f^{-1}W$ is a union of multiplicity-strata
of $f^{-1}(V\cap M)$, each of which is made by~$f$ into a $C^\infty\<$-fiber
bundle, with smooth fibers, over~$W$? We could further require that the
not-necessarily-reduced space~$f^{-1}W$ be locally analytically trivial
over~$W$.

As mentioned above, Zariski gave an affirmative answer when $\dt x=1$.
We outline now an approach to the question which most likely gives an
affirmative answer when $\dt x=2$.
(At this writing, I have checked many, but not all, of the details.)
Roughly speaking, this approach is the classical desingularization method of
Jung, applied to families.  The main roadblock to extending it inductively to
dimensionality types $\ge 3$ is indicated  below (Problem~\ref{QOres}).

For simplicity, assume that $\dim V=3$ and that $\dt x=2$, so that the Zariski
2-stratum is a non-singular curve~$W$ through~$x$. By definition
we can choose a finite map germ  whose discriminant $\Delta$ has the image
of~$W$ as its Zariski 1-stratum; and after reimbedding $V$ we may
assume this map germ to be induced by a linear map $\pi\colon\C^{\>4}\to
\C^{\>3}$. In what follows, we need only \emph{one} such~$\pi$.

According to Zariski, there is an embedded resolution $h\colon M^3\to\C^{\>3}$
of~$\Delta$ such that $h^{-1}\pi(W)$ is a union of strata.
(We abuse notation by writing $\C^{\>3}$ for a suitably small
neighborhood of~$\pi(x)$ in~$\C^{\>3}$.) Let
$g\colon M^4 = M^3\times_{\C^{3}} \C^{\>4}\to \C^{\>4}$ be the projection, a
birational map of manifolds. There is a finite map
from the codimension-1 subvariety $V'\!:=g^{-1}V\subset M^4$ to
$M^3$ whose discriminant, being contained in $h^{-1}\Delta$, has normal
crossings. This means that the singularities of~$V'$ are all {\it
quasi-ordinary}, see \cite{L2}. Moreover, as Zariski showed (see
\cite[p.\,514]{Z1}), the Zariski-strata of~$V$ are all \'etale over the
corresponding strata of~$\Delta$, so their behavior under pullback through
$V'\to V$ is essentially the same as that of the strata of $\Delta$ under~$h$,
giving rise to nice fiber-bundle structures. More precisely, for the
Zariski-stratification~$\Sigma$ on~$V$ and for the pullback~$\Sigma'$ on
$V'\<$ of the multiplicity stratification on the normal crossings divisor
$g^{-1}\Delta$,  $g$ makes each
$\Sigma'$-stratum into a $C^{\infty\<}$-fiber bundle over a $\Sigma$-stratum.
(In other words, the map~$g$ is \emph{stratified} with respect to the
indicated stratifications.) Thus we have achieved an ``equisimplification," but
not
yet an equiresolution---the singularities of~$V$ along~$W$ have simultaneously
been made quasi-ordinary.

So now we need only deal with quasi-ordinary singularities,
keeping in mind however that $V'$ is no longer a localized object---it
contains projective subvarieties in the fibers over the original
singularities of~$V\<$. Fortunately, the above-defined
stratification $\Sigma'$ can be characterized \emph{intrinsically} (even
\emph{topologically}) on~$V'\<$, see \cite [\S6.5]{L3}.  So we can forget about
the
projection $V'\to M^3$ via which $\Sigma'$ was determined, and deal directly
with this canonical stratification on~$V'$ (closely related, presumably, with
the
Zariski stratification, though
$V'\to M^3$ may not be sufficiently general at every point of~$V'$).  The
problem is thus reduced to showing that \emph{the canonical
stratification\/~$\Sigma'$ on  the quasi-ordinary space\/~$V'$ is
equiresolvable.}

\smallskip
Any germ $(V'\<,y)$ (assumed, for simplicity, irreducible) is given by the
vanishing
of a polynomial
$$
Z^m+a_1(X,Y,t)Z^{m-1}+\cdots+a_m(X,Y,t)\qquad (a_i\in\C[[X,Y,t]])
$$
whose discriminant is of the form
$$
\delta=X^aY^b\e(X,Y,t),\qquad \e(0,0,0)\ne 0.
$$
(Here $t$ is a local parameter along~$W$.) The roots of this
polynomial are fractional power series $\zeta_i(X^{1/n}\<, Y^{1/n}\<,t)$;
and since $\delta=\prod(\zeta_i-\zeta_j)$, we have, for some
non-negative integers $a_{ij}, b_{ij}$,
$$
\zeta_i-\zeta_j=X^{a_{ij}/n}Y^{b_{ij}/n}\e_{ij}(X^{1/n}\<,Y^{1/n}\<,t),
\qquad \e_{ij}(0,0,0)\ne 0.
$$
Modulo some standardization, the monomials $X^{a_{ij}/n}Y^{b_{ij}/n}$ so
obtained are called the {\it characteristic monomials\/} of the quasi-ordinary
germ~$(V'\<,y)$. They provide a very effective tool for studying such germs.
They are higher dimensional generalizations of the characteristic
pairs of  plane curve singularities (which are always quasi-ordinary),
and they control many basic features of $(V',y)$, for example the
number of components of the singular locus, the multiplicities of
these components on $V'$ and at the origin, etc. (See \cite{L2}, \cite{L3} for
more details.) In particular, two quasi-ordinary singularities have the
same embedded topology iff they have the same characteristic monomials
\cite{G}.
It is therefore natural to ask:


\begin{problem} \label{QOres}

Do the characteristic monomials of a
quasi-ordinary singularity determine a canonical embedded resolution
procedure?

\end{problem}


This question makes sense in any dimension, and would arise naturally in
any attempt to extend the preceding argument inductively to higher
dimensionality types.  A positive answer would imply that the
canonical stratification of any quasi-ordinary hypersurface is equiresolvable.

At this point,
if we weren't concerned with embedded resolution (which we need to be, if
there is to be any possibility of induction), we could  \emph{normalize}~$V'$
to get a family of cyclic quotient singularities whose structure is completely
determined by the characteristic monomials. Simultaneously resolving such a
family is old hat (cf.~\cite[Lecture 2]{L1}),   leading to the following result
(for germs of families of two-dimensional hypersurface singularities): \emph{If
there is a single projection with an equisingular branch locus, then we have
simultaneous resolution; and if that projection is transversal (i.e., its
direction
does not lie in the Zariski tangent cone) then we have
strong simultaneous resolution, and in
particular (by Teissier's result) differential equisingularity. }

In particular, the
above-mentioned example of Luengo \emph{is} in fact strongly simultaneously
resolvable---though not by one of the standard blowup methods.

For \emph{embedded resolution,} canonical algorithms are now available
\cite{BM}, \cite{EV}, so that Problem~\ref{QOres} is quite concrete: verify
that the invariants which drive a canonical procedure, operating on a
quasi-ordinary singularity parametrized as above (by the~$\zeta_i$), are
completely determined by the characteristic monomials.

\penalty -2000

The idea is then that the intersections of $V'$ with germs of smooth
varieties transversal to the strata form families of quasi-ordinary
singularities with the same characteristic monomials, and so any
resolution of one member of the family should propagate along the
entire stratum---whence the equiresolvability

A key point in dimension~2 is {\it monoidal stability:} any permissible
blowup of a 2-dimensional quasi-ordinary\ singularity is
again quasi-ordinary, and the
characteristic monomials of the blown up singularity depend only on
those of the original one and the center of blowing up---chosen according to an
algorithm depending only on the characteristic monomials; explicit formulas are
given in \cite[p.\,170]{L2}.  When it comes to total transform, there are of
course
some complications; but monoidal stability can still be worked out for a
configuration consisting of a quasi-ordinary singularity \emph{together with} a
normal crossings divisor such as arises in the course of embedded resolution.
In principle, then, one
should be able to use an available canonical resolution process and see things
through. Some preliminary work along these lines was reported on at the
workshop by Chungsheng Ban and Lee McEwan.

Unfortunately such monoidal stability fails in higher dimensions, even for
such simple quasi-ordinary singularities as the origin on the threefold
$W^4=XYZ$.  So we need to
look into:


\begin{problem}

Find a condition on singularities weaker than
quasi-ordinariness, but which is monoidally stable, and which can be
substituted for quasi-ordinariness in~Problem~\ref{QOres}.

\end{problem}


Careful analysis of how the above-mentioned canonical resolution procedures
work on  quasi-ordinary singularities could suggest an answer.



\bigskip\bigskip\bigskip


\noindent
Department of Mathematics\newline Purdue University\newline
W. Lafayette, IN 47907, USA

\medskip

\noindent
lipman@math.purdue.edu\newline
www.math.purdue.edu/\~{}lipman

\end{document}